# Memory of Ivan Kupka (1937-2023)

## D. Holcman

**Abstract:**

This brief text is in memory of Professor Ivan Kupka. It presents his vision, scientific life, his interest in mathematics and our join collaboration.

## Trajectoire furieuse d'une vocation mathématique

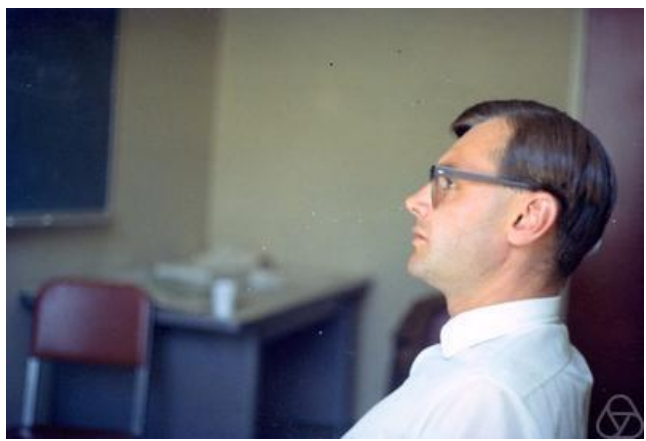

Location: Berkeley
Author: Bergman, George M. (photos provided by Bergman, George M.)
Source: George M. Bergman, Berkeley
Year: 1968
Copyright: George M. Bergman, Berkeley

Voici l'histoire brève d'un homme qui a traversé la vie comme un cyclone dans l'univers des mathématiques, mais aussi à travers les pays, sans s'arrêter, sans frontières avec ces 6 passeports, en oscillant entre les contraintes. Ni besoin de gloire ou de reconnaissance, animé jusqu'à la fin dans son Epadh de sèvre où il avait perdu la mémoire, du besoin de comprendre, de calculer, de réfléchir mais il s'amusait. Un homme qui était généreux de son temps avec ceux qui s'intéressaient aux math, à la physique, à la chimie, et à la science en général.

Kupka, il se présentait comme ça, métèque me disait-il, quand il parlait de lui-même. Il est né en 1937, dans la ville de Třebíč en Tchécoslovaquie. Bien que la famille séjourne en France avant 1939, ils sont retournés en Tchécoslovaquie pour survivre, et passer la Seconde Guerre mondiale. La famille est revenue en France après, pour échapper au communisme d'après-guerre. Ivan est entré à l'ENS sans faire de math sup, après seulement une année de préparation. Il racontait comment il expliquait à ses camarades surtout, de lettres de l'ENS, ce que le communisme faisait à la Techchoslovique, bien que cela n'arrivât pas à briser l'endoctrinement de plusieurs générations. Il est ensuite parti pour voir autre chose, mais il aimait bien l'ENS. Alors, il traversa l'Amérique du Sud pendant 2 ans, avec son calibre 22, qui lui servira pour survivre en Colombie à l'époque du sentier lumineux. Il fera des petits métiers, pécheurs, réparations des moteurs de bateaux pour la pêche à la crevette à l'est du Pérou. Ces années d'errance lui ont permis de confirmer sa mission pour les maths.

Un jour, dans un de ces diners de con entre X-promo-Années, comme ils se présentaient bien, assis à coté de Jean-Louis Basdevant, physicien émérite de l'école Polytechnique, je lui demande s'il connaît Kupka. Il me répondit: bien sûr, impossible d'oublier celui qu'il redécouvrit à Berkeley en tant que cuistot sur la liaison maritime Lima-San Francisco. La nouvelle avait circulé très rapidement à Berkeley, qu'un cuisinier savait résoudre n'importe quel problème de math. Finalement, après avoir été reconnu, M. Peixoto, au Brésil allait lui obtenir une bourse de thèse. Le voilà pour le brésil, où il découvrira la femme dans toute ça splendeur. Au même moment, pendant sa thèse, les autorités françaises recherchaient Ivan pour qu'il rejoigne l'armée française engagée en Algérie. Ne s'étant finalement pas présenté au service militaire, il restera au Brésil pour finir sa thèse et sera considéré comme insoumis, statut qui lui interdira de revenir en France pendant 22 ans.

Dès ses débuts, Ivan aura des contributions significatives dans les systèmes dynamiques [Théorème de Kupka-Smale sur la densité des champs de vecteurs en dimension 2, dont les singularités sont hyperboliques]. Ces contributions se retrouvent en géométrie riemannienne et sous-Riemannienne, en théorie du contrôle, Équations aux dérivées partielles, méthodes numériques, symplectique, et physique mathématiques pour ne citer que les principales. Ivan occupa aussi diverses positions universitaires à travers les grandes universités prestigieuses: de sa thèse au Brésil sous la direction de Peixoto, il fut ensuite professeur à Berkeley, Stony-Brook, Toronto, Grenoble et Paris. Chern n'aimait pas Ivan, et Ivan ne correspondait pas bien dans la tradition chinoise du successeur du maitre. Donc pas de « tenure » à Berkeley. Il ira donc à Stony Brook, côte Est, ou Jim Simons était Chairman du département.

J'ai rencontré Ivan lorsque j'étais étudiant à l'Université Paris VI, lors de son cours sur la géométrie symplectique en 1996-97. En quelques années, nous sommes devenus amis et avons commencé à travailler ensemble sur tous les sujets qui nous intéressaient. Ivan pouvait discuter librement de tous les points de mathématiques, il avait une vision encyclopédique et dévouée, rare. Pas de conneries avec Ivan, pas de superflus, pas de chichi, l'essentiel, les discussions étant dans le sujet. Il allait droit au but, résumant en une phrase les notions mathématiques abstraites de toutes sortes. Ivan était un mathématicien hors norme, qui communiquait sa passion des mathématiques avec énergie, directe. Combien de calculs dans tous les sens, rien ne l'arrêtait, il était spécialiste en toutes ces sous-disciplines. Une encyclopédie au service du travail de recherche; Je n'ai pas souvenir d'une théorie mathématique qui ne l'avait pas intéressé. Un esprit curieux, aventurier, rebelle et dédier. Ivan travaillait sans cesse, le jour, le soir le week-end, les vacances. C'était un plaisir et privilège de travailler avec lui.

Nous avons travaillé 20 ans ensemble sans s'arrêter, toujours animés par l'envie de comprendre et de calculer. Nous avons commencé à travailler sur les EDP du premier ordre sur les variétés Riemanniennes, où nous avons découvert que les solutions dépendaient des points critiques ou des cycles limites associés aux systèmes dynamiques (terme du premier ordre), ce qui n'intéressa personne. Nous avons poursuivi avec les limites des perturbations singulières de la première valeur propre du Laplacien sur les variétés riemanniennes, où la limite était liée à la pression topologique, une notion développée par Y. Kiefer dans les années 80-90. Nous avons ensuite changé complètement de direction. Explorant le problème de deux particules browniennes sur un intervalle: quelle était la probabilité que l'une d'entre elles s'échappe avant qu'elle se heurtent ? la probabilité était liée à la fonction de Weierstrass. Un résultat qui avait intrigué l'insatiable Marc Yor [décédé en 2013] dans les calculs des probabilités. Nous avons ensuite étudié le premier temps de rencontre de 2 extrémités d'un polymère de Rouse (N perles reliées par des ressorts). Ce problème nous a conduit à utiliser la théorie de Chavel-Fledman de la perturbation pour la première valeur propre du Laplacien en haute dimension sur un sous-variété avec un petit voisinage tubulaire. Dans les années 2010-12, nous avons trouvé un développement asymptotique de la solution par rapport à la dimension, qui était ici le nombre de billes du polymère. Nous avons obtenu des estimations asymptotiques pour le temps de formation d'une boucle, ce qui était pour nous une grande satisfaction. Il y a eu plusieurs tentatives d'obtenir des formules asymptotiques par Pastor au NIH dans les années 90. Ces formules sont devenues utiles avec le développement du domaine de l'organisation de l'ADN dans le noyau cellulaire, domaine qui a émergé à l'interface entre la physique statistique, la théorie des polymères en grande dimension, les processus stochastiques, l'asymptotique mathématique et la génomique tri-dimensionnelle.

Nous avons également étudié le temps nécessaire aux spermatozoïdes pour trouver un ovule dans l'utérus. Nous avons commencé par un modèle de type billard, où la trajectoire est linéaire à l'intérieur d'un domaine, mais la réflexion sur la frontière est aléatoire : combien de temps faut-il pour trouver un petit bout du bord ? Nous avons également étudié le temps du plus rapide. La trajectoire du plus rapide est en fait donnée par un problème de minimisation. La trajectoire du plus rapide est donc très différente de celle d'un spermatozoïde qui arrivera bien plus tard dans le temps. La trajectoire du plus rapide est sélectionnée par la statistique extrême. Ainsi, les propriétés de la trajectoire physique initiale, pour autant qu'elle soit suffisamment aléatoire, ne sont pas nécessairement pertinentes pour le premier temps d'arrivé (à un ovule). Cela a ouvert la question générale de savoir pourquoi il y a tant de spermatozoïdes et qu'est-ce qui en détermine ce grand nombre dans la nature. Notre réponse était qu'il faut compenser pour trouver une petite portion du domaine où se trouve l'œuf à fertiliser. Cela à donner naissance a la théorie de la redondance en biologie.

Pendant toutes ces années, nous avons eu beaucoup de plaisir à travailler ensemble, à trouver des problèmes, à discuter des méthodes de calcul. Cela m'a permis de définir les valeurs mathématiques à communiquer. Ivan n'était pas intéressé par les fonctions administratives de petit et haut niveau, il était simplement intéressé par les mathématiques, mais il nous a poussés à publier, à écrire des demandes de financements pour survivre dans ce XXI siècle, rempli de contraintes administratives en tous genres que l'informatique aurait dû faire disparaitre.

Ivan a travaillé autant qu'il a pu, mais il a souffert, au cours des deux dernières années, de pertes de mémoire profondes, en particulier après le covid. Il est décédé en avril 2023. La cérémonie fut courte avant de passer au crématorium. Nous étions au maximum 10 à ses funérailles, la famille et trois collègues. Il a fallu aussi abréger car le suivant attendait, on a donc réduit les trois discours. Cette cérémonie, je l'avais déjà vu il y 35 ans, c'était la même que la scène d'Amadeus tombant dans la fosse publique https://www.youtube.com/watch?v=TUt4DfGnyJQ C'était donc un classique que d'en finir rapidement et sans spectateurs pour ceux qui ont travaillé dure. Pourtant cette fin-là, c'est aussi une génération de passionnés qui s'en va, montrant que la vie des scientifiques a bien changé. Alors qui va se charger des grandes questions ? l'IA ? ou la nouvelle génération d'étudiants dont l'idéal est maintenant le Machine-learning et tiquer sur leur téléphone. Le temps de la recherche s'est réduit, dû aux exercices administratifs, plus de secrétaires, plus de gestionnaires, à quoi bon dans ce grand déguisement. Kupka, c'était rien de tous ca. Sa façon de dire : travailler avec passion, c'est la liberté que nous devons avoir, un style de vie, pour ces questions nouvelles, sans besoin de rien d'autres, qui a l'avant-garde permet de faire avancer la connaissance. Comme on le faisait avec Ivan, il faut chaque jour peaufiner sa technique de calcul, les idées, les modèles, les théories, les simulations, rien négliger.

Kupka reste pour beaucoup d'entre nous un modèle original portant des valeurs du savoir, de liberté intellectuel, de travail et d'amitiés.

Je voudrais terminer avec cette anecdote de Charles Pugh

*************************************************************************

I have fond memories of Ivan's eccentricities. I'll give you one from about 1965. My wife and I invited Ivan to lunch. He was three hours late. No apologies. We were glad to see him. We were sitting around the table chatting, and I was admiring his glasses — the French frameless style was unknown in the US at the time. Although stylish and wonderfully minimal, I expressed doubts about their physical strength. After some theoretical to and froh about the weakness of the wire connection to the glass, and without any pause in the conversation, Ivan flung them at the wall across the room. No damage. No glee. Case closed. I loved his boredom with unfounded opinion.

*Charles Pugh*, UC Berkeley.
. . . . . . . . . . . . . . . . . . . . . . . . . . . . . . . . . . . . . . . . . . . . . . . . . . . . . . . . . . . . . . . . . . . . . . . . . . . . . . . . . .

Et cet Hommage de Dennis Sullivan

I met Ivan Kupka during the late sixties and early seventies while becoming interested in dynamics after working in algebraic topology. His beautiful result with Smale that for a generic dynamical system, an orbit is either wandering or approximated by hyperbolic periodic orbits whose attractive and expanding leaves all intersect transversally gave a vivid picture to build on for decades . Ivan ,whenever we met, was very helpful in this process.

Dennis Sullivan
*************************************************************************